\documentclass[11pt]{article}

\usepackage{amsfonts, amsthm, amssymb}

\theoremstyle{plain}

\newtheorem*{remark}{Remark}
\newtheorem*{example}{Example}
\newtheorem{thm}{Theorem}
\newtheorem{prop}{Proposition}[section]
\newtheorem{cor}[prop]{Corollary}

\newtheorem{step}{Step}

\newcounter{ListingNr}

\newenvironment{listing}%
{\begin{list}   { {\rm (\roman{ListingNr})} }
                { \topsep2mm
                  \partopsep0ex
                  \listparindent0pt
                  \itemsep0.5ex
                  \parskip0pt
                  \leftmargin5ex
                  \labelwidth5.5ex
                  \labelsep0ex
                  
                  \usecounter{ListingNr}} }%
{\end{list}}

\begin{document}

%
%

\newcommand{\emptyword}{\varnothing}
\newcommand{\length}[1]{\ell(#1)}
\newcommand{\shuffle}{\unitlength1pt%
  \begin{picture}(10,8)%
    \linethickness{0.4pt}
    \put(2,0){\line(1,0){6}}%
    \put(2,0){\line(0,1){4}}%
    \put(5,0){\line(0,1){4}}%
    \put(8,0){\line(0,1){4}}%
  \end{picture}}

\newcommand{\relLK}{\equiv_{\mathsf{LQ}}}
\newcommand{\relL}{\equiv_{\mathsf{L}}}
\newcommand{\relK}{\equiv_{\mathsf{Q}}}

\newcommand{\N}{{\mathbb N}}
\newcommand{\Z}{\mathbb{Z}}
\newcommand{\Q}{\mathbb{Q}}

\newcommand{\Sym}{\mathcal{S}}

\newcommand{\CC}{\mathcal{C}}
\newcommand{\FF}{\mathcal{F}}
\renewcommand{\AA}{\mathcal{A}}
\newcommand{\QQ}{\mathcal{Q}}
\newcommand{\DD}{\mathcal{D}}
\newcommand{\LL}{\mathcal{L}}
\newcommand{\PP}{\mathcal{P}}

\newcommand{\Lie}{\mathsf{Lie}}
\newcommand{\Des}{\mathsf{Des}}
\newcommand{\Prim}{\mathsf{Prim}}

\newcommand{\setofall}[2]{\mbox{$\{\,#1\,|\,#2\,\}$}}
\newcommand{\Setofall}[2]{\mbox{$\Big\{\,#1\,\Big|\,#2\,\Big\}$}}
\newcommand{\ol}{\overline}
\newcommand{\ul}{\underline}
\newcommand{\Erz}[1]{\Big\langle\,#1\,\Big\rangle}
\def\haken#1{\underline{#1}%
             {\raise -0.35ex%
              \hbox{\vphantom{$#1$}\vrule height 0.8ex}}}
\newcommand{\kref}[1]{(\ref{#1})}
\newcommand{\tensor}{\otimes}

%
%

%
%

\title{
Lie Elements and Knuth Relations}

\author{
Manfred Schocker\thanks{supported by the Research Chairs of Canada}\\[1mm]
  \begin{minipage}{5.5cm}
    \begin{center}
         \begin{small}
           LaCIM, UQAM\\
           CP 8888, succ. Centre-Ville\\
           Montr{\'e}al (Qu{\'e}bec), H3C 3P8\\
           Canada\\
           E-mail: \texttt{mschock@math.uqam.ca}
         \end{small}
    \end{center}
  \end{minipage}
}

\maketitle

%
%

\begin{abstract} \noindent
  A coplactic class in the symmetric group $\Sym_n$ consists of all
  permutations in $\Sym_n$ with a given Schensted $Q$-symbol, and may
  be described in terms of local relations introduced by Knuth. Any
  Lie element in the group algebra of $\Sym_n$ which is constant on
  coplactic classes is already constant on descent classes. As a
  consequence, the intersection of the Lie convolution algebra
  introduced by Patras and Reutenauer and the coplactic algebra
  introduced by Poirier and Reutenauer is the Solomon descent algebra.
\end{abstract}

\begin{center}
\end{center}

%
%

\section{Introduction}

In 1995, Malvenuto and Reutenauer introduced the structure of a graded
Hopf algebra on the direct sum
$$
\PP
=
\bigoplus_{n\ge 0} \Z\Sym_n
$$
of all symmetric group algebras $\Z\Sym_n$ over the ring $\Z$ of
integers (\cite{malvenuto-reutenauer95}).  Apart from this
\emph{convolution algebra of permutations} $\PP$ itself
(\cite{aguiar-sottile02, parisVI}) several sub-algebras of $\PP$
turned out to be of particular algebraic and combinatorial interest
and have been studied intensively; for instance, the Rahmen\-algebra
(\cite{joellenbeck99}), the Hopf algebra of the planar binary trees
(\cite{loday-ronco98,chapoton00}), the Lie convolution algebra $\LL$
(\cite{patrasreut01}), the coplactic algebra $\QQ$
(\cite{poirierreute95}\footnote{%
  The algebra $(\mathbb{Z}C,*,\delta)$ introduced in \cite{poirierreute95}
  is the dual algebra of the algebra $\QQ$ considered here
  (see \cite[Th\'eor\`eme 3.4]{poirierreute95}).},
\cite{blescho02}), and the Solomon descent algebra $\DD$
(\cite{%
  solomon76,garsia-reutenauer89,reutenauer93,malvenuto-reutenauer95,%
  parisI,blelau96,joellereut01}).  

Here, the \emph{relation} between the algebras $\LL$ and $\QQ$ shall
be investigated.  The latter is defined combinatorially as the linear
span of the sums of permutations with given Schensted $Q$-symbol
(\cite{schensted61}), or, equivalently, of the sums of equivalence
classes arising from the coplactic relations in $\Sym_n$, $n\ge 0$,
introduced by Knuth (\cite{knuth70}).  The Lie convolution algebra
$\LL$ is generated (as an algebra) by all Lie elements in $\PP$. Both
$\LL$ and $\QQ$ contain $\DD$. Combinatorial descriptions of the
algebras $\DD$ and $\QQ$, and the set of Lie elements in $\PP$, follow
in Section~\ref{relations}.  The main goal of this paper is to show

\begin{thm} \label{intersect}
  $\LL\cap \QQ=\DD$.
\end{thm}

This result (once more) points out the exceptional role played by the
Solomon algebra. The proof is given in Section~\ref{relations}, and is
essentially based on the fact that any Lie element in $\PP$ which is
constant on coplactic classes is already contained in $\DD$ (see
Section~\ref{copl-lie}), which is combinatorially interesting for its
own sake.

One might be tempted to conjecture that a lack of co-commutativity of
$\QQ$ is the deeper reason for Theorem~\ref{intersect}, since $\LL$ is
--- at least in comparison to $\DD$ --- a ``large'' co-commutative
sub-algebra of $\PP$; but this is false. The domain of co-commutativity
of $\QQ$ \emph{strictly} contains $\DD$.  Some comments concerning
this can be found at the end of Section~\ref{copl-lie}.

%
%

\section{Descent, Coplactic, and Lie Relations} \label{relations}

In this section, combinatorial descriptions of the algebras $\DD$ and
$\QQ$, and of the Lie elements in $\PP$, are recalled briefly; and a
proof of Theorem~\ref{intersect} is given.

Let $\N$ (respectively, $\N_0$) be the set of positive (respectively,
nonnegative) integers and set
$$
\haken{n}
:=
\setofall{i\in\N}{i\le n}
$$
for all integers $n$.  For any $\pi\in \Sym_n$,
$
\Des(\pi)
:=
\setofall{i\in\haken{n-1}}{i\pi>(i+1)\pi}
$
is the \emph{descent set} of $\pi$.  The Solomon descent algebra
$\DD$ is the linear span of the sums
$
\sum_{\pi\in \Sym_n\atop \Des(\pi)=D} \pi,
$
where $n\in\N_0$ and $D\subseteq \haken{n-1}$.  Due to Malvenuto
and Reutenauer, $\DD$ is a Hopf sub-algebra of $\PP$
(\cite{malvenuto-reutenauer95}).  We mention that
the homogeneous component $\DD_n=\DD\cap \Z\Sym_n$ is a sub-algebra of
$\Z\Sym_n$, due to a remarkable result of Solomon (\cite{solomon76}),
although this is not of relevance here.

Let $\N^*$ be a free monoid over the alphabet $\N$ and denote by
$\emptyword$ the empty word in $\N^*$.  The mapping $\pi\mapsto
(1\pi)\ldots(n\pi)$ extends to a linear embedding of $\Z\Sym_n$ into the
semi-group algebra $\Z\N^*$.  As is convenient for our purposes,
elements of $\Z\Sym_n$ will be identified with the corresponding elements
of $\Z\N^*$.  Furthermore, products $\sigma\nu$ of permutations
$\sigma,\nu\in \Sym_n$ are to be read from left to right: first $\sigma$,
then $\pi$.

The following combinatorial characterization of $\DD_n$
was given in \cite[4.2]{blelau93a}.

\begin{prop}[Descent Relations] \label{descent-char}
  Let $\varphi=\sum_{\nu\in \Sym_n} k_\nu \nu\in \Z\Sym_n$, then
  $\varphi\in\DD_n$ if and only if
  $$
  k_{uaw(a+1)v} = k_{u(a+1)wav}
  $$
  for all $a\in\haken{n-1}$, $u,v,w\in\N^*$ such that $\nu=uaw(a+1)v\in
  \Sym_n$ and $w\neq\emptyword$.
\end{prop}

Let $Q(\pi)$ denote the Schensted $Q$-symbol of $\pi$,
for all $\pi\in\Sym_n$ (\cite{schensted61}), then the set
of all $\sigma\in\Sym_n$  such that 
$Q(\pi)=Q(\sigma)$ is a \emph{coplactic class} in
$\Sym_n$\footnote{%
  Due to Sch\"utzenberger (\cite{schuetzen63}), $P(\pi)=Q(\pi^{-1})$
  is the Schensted $P$-symbol of $\pi$; and the equivalence arising
  from equality of $P$-symbols leads to the \emph{plactic monoid}
  (\cite{lasschue81}).  This is the reason why the word coplactic is
  used here.}.  The coplactic algebra $\QQ$ is the linear span of all
sums of coplactic classes in $\PP$:
$$
\QQ
=
\Erz{%
  \Setofall{\sum_{Q(\sigma)=Q(\pi)}\sigma}{\pi\in\Sym_n,\; n\in\N_0}}_\Z\,.
$$
Accordingly, each element $\varphi\in\QQ$ is called coplactic.  Due
to Poirier and Reutenauer, $\QQ$ is a Hopf sub-algebra of $\PP$
(\cite{poirierreute95}).  The following characterization of
$\QQ_n:=\QQ\cap\Z\Sym_n$ is due to Knuth (\cite{knuth70}).

\begin{prop}[Coplactic Relations] \label{coplactic-char}
  Let $\varphi=\sum_{\nu\in \Sym_n} k_\nu \nu\in \Z\Sym_n$, then
  $\varphi\in\QQ_n$ if and only if
  $$
  k_{uaw(a+1)v} = k_{u(a+1)wav}
  $$
  for all $a\in\haken{n-1}$, $u,v,w\in\N^*$ such that $\nu=uaw(a+1)v\in
  \Sym_n$ and $w$ contains the letter $a-1$ or the letter $a+2$.
\end{prop}

Combining Propositions \ref{descent-char} and \ref{coplactic-char}
implies, in particular, $\DD\subseteq \QQ$.

Let
$$
\omega_n
=
\sum_\pi (-1)^{1\pi^{-1}-1} \pi
\in \Z\Sym_n\,,
$$
where the sum is taken over all \emph{valley permutations} $\pi\in
\Sym_n$, which are defined by the property
$1\pi>\cdots>(k-1)\pi>k\pi<(k+1)\pi<\cdots<n\pi$, where
$k:=1\pi^{-1}$.  The element $\omega_n$ projects $\Z\Sym_n$ onto the
multi-linear part of the free Lie algebra, by left multiplication
(\cite{dynkin47,specht48,wever49}, see \cite{blelau93a}). Accordingly,
$$
\Lie_n
:=
\omega_n\Z\Sym_n
$$
is the set of \emph{Lie elements} in $\Z\Sym_n$ for all $n\in\N_0$.
Each 
$\varphi\in\Lie:=\bigoplus_{n\ge 0} \Lie_n$
is a primitive element of the Hopf algebra $\PP$ (\cite{patrasreut01}).
The Lie convolution algebra $\LL$ is the (co-commutative) Hopf sub-algebra
of $\PP$ generated by $\Lie$; there is also the relation
$\DD\subseteq\LL$ (\cite{patrasreut01}).

In view of a \textbf{proof of Theorem~\ref{intersect}}, consider the
corresponding algebras $\DD_\Q$, $\LL_\Q$, $\QQ_\Q$, and $\PP_\Q$ over
the field $\Q$ of rational numbers, then $\DD_\Q$ is contained in
$\LL_\Q\cap\QQ_\Q$; the latter is a co-commutative Hopf sub-algebra of
$\PP_\Q$, hence generated by its primitive elements, due to Milnor and
Moore (\cite{milnor-moore65}).  But each primitive element in
$\LL_\Q\cap\QQ_\Q$ is, in particular, a primitive element in $\LL_\Q$
and therefore (up to a rational factor) contained in $\Lie$.  In
Section~\ref{copl-lie}, it will be shown that any coplactic Lie
element $\varphi\in\Lie\cap\QQ$ is contained in $\DD$
(Theorem~\ref{copl-lie-thm}).  This implies
$\LL_\Q\cap\QQ_\Q\subseteq\DD_\Q$.  Observing that
$\DD_\Q\cap\PP=\DD$, completes the proof of Theorem~\ref{intersect}.

A combinatorial characterization of the set $\Lie_n$ follows.  Let
$u\shuffle v$ denote the usual shuffle product of $u=u_1\ldots
u_k,v=v_1\ldots v_m\in\N^*$, that is
$$
u\shuffle v=\sum_w w,
$$
where the sum ranges over all $w=w_1\ldots w_{m+k}\in\N^*$ such
that $u=w_{i_1}\ldots w_{i_k}$ and $v=w_{j_1}\ldots w_{j_m}$ for
suitably chosen indices $i_1<\cdots<i_k$, $j_1<\cdots<j_m$ such that
$\haken{k+m}=\{i_1,\ldots,i_k,j_1,\ldots,j_m\}$.
Furthermore, set
$$
\ol{u}
:=
u_k\ldots u_1
$$
and denote by $\length{u}:=k$ the length of $u$.

\begin{prop} \label{prepare}
  Let $n\in\N$ and $a\in\haken{n}$, then 
  $\setofall{\omega_n\sigma}{\sigma\in \Sym_n,\,1\sigma=a}$
  is a linear basis of $\Lie_n$. 
  
  Furthermore, for any choice of coefficients $c_\sigma\in \Z$
  ($\sigma\in \Sym_n$, $1\sigma=a$), the coefficient of $\nu=uav\in \Sym_n$ in
  $\omega_n\sum_{1\sigma=a}c_\sigma\sigma$ is
  \begin{equation}
    \label{shuffle-formula}
    (-1)^{\length{u}}
    c_{a(\ol{u}\shuffle v)}\,,    
  \end{equation}
  where $\pi\mapsto c_\pi$ has been extended to $\Z\Sym_n$ linearly.  In
  particular, the coefficient of $\sigma\in \Sym_n$ 
  is $c_\sigma$ whenever $1\sigma=a$.
\end{prop}

This result is seemingly folklore; a proof follows for the reader's
convenience.

\begin{proof}
  Let $\nu=uav\in \Sym_n$ and $\sigma=ax_2\ldots x_n\in \Sym_n$, then the
  coefficient of $\nu$ in $\omega_n\sigma$ is non-zero if and only if
  there is a valley permutation $\pi\in \Sym_n$ such that
  $$
  uav
  =
  \nu
  =
  \pi\sigma
  =
  x_{1\pi}\ldots x_{(k-1)\pi}ax_{(k+1)\pi}\ldots x_{n\pi}\,,
  $$
  where $k:=1\pi^{-1}$; that is,
  $u=x_{1\pi}\ldots x_{(k-1)\pi}$ 
  and
  $v=x_{(k+1)\pi}\ldots x_{n\pi}$.
  Since $1\pi>\cdots>(k-1)\pi$ and $(k+1)\pi<\cdots<n\pi$,
  this is equivalent to saying that $x_2\ldots x_n$ is a summand in
  the shuffle product of $\ol{u}$ and $v$;
  in this case, the coefficient of $\nu$ in $\omega_n\sigma$ is
  $(-1)^{1\pi^{-1}-1}=(-1)^{\length{u}}$.
  This proves \kref{shuffle-formula}. Since 
  $$
  \dim\,\Lie_n
  = 
  (n-1)! 
  = \#\setofall{\omega_n\sigma}{\sigma\in \Sym_n,\,1\sigma=a}
  $$
  and the coefficient of $\tilde{\sigma}=av\in \Sym_n$ in
  $\omega_n\sum_{1\sigma=a}c_\sigma\sigma$ is $c_{\tilde{\sigma}}$,
  the basis property follows.
\end{proof}

\begin{cor}[Lie relations] \label{lie-char}
  Let $\varphi=\sum_{\nu\in \Sym_n} k_\nu \nu\in \Z\Sym_n$, then
  $\varphi\in\Lie_n$ if and only if
  \begin{equation}
    \label{lie-shuffle}
    k_{uav} = (-1)^{\length{u}} k_{a(\ol{u}\shuffle v)}
  \end{equation}
  for all $a\in\haken{n}$, $u,v\in\N^*$ such that $\nu=uav\in \Sym_n$.
\end{cor}

\begin{proof}
  Let $\varphi\in\Lie_n$ and $a\in\haken{n}$, then there are 
  coefficients $c_\sigma\in \Z$
  ($\sigma\in \Sym_n$, $1\sigma=a$) such that
  $\varphi=\omega_n\sum_{1\sigma=a}c_\sigma\sigma$,
  by Proposition~\ref{prepare}, and
  $$
  k_{uav}
  =
  (-1)^{\length{u}}
  c_{a(\ol{u}\shuffle v)}
  =
  (-1)^{\length{u}}
  k_{a(\ol{u}\shuffle v)}\,,
  $$
  by \kref{shuffle-formula}.
  Conversely, \kref{lie-shuffle} implies
  $
  \varphi
  =
  \omega_n\sum_{1\sigma=a}k_\sigma\sigma \in\Lie_n
  $,
  by \kref{shuffle-formula} again.
\end{proof}

Proposition~\ref{coplactic-char} and Corollary~\ref{lie-char} may be
restated as follows.  Let $T^\perp$ be the space orthogonal to $T$
with respect to the scalar product on $\Z\Sym_n$ turning $\Sym_n$ into
an orthonormal basis, for all $T\subseteq\Z\Sym_n$. For all
$c,d\in\Z\Sym_n$, write
\begin{center}
$c\relK d$ (respectively, $c\relL d$, $c\relLK d$),
\end{center}
if
$c-d\in \QQ_n^\perp$ (respectively,
$\in \Lie_n^\perp$, $\in (\Lie_n\cap\QQ_n)^\perp$).
Now the necessity parts of Proposition~\ref{coplactic-char}
and Corollary~\ref{lie-char} are
\begin{equation}
  \label{relq}
  uaw(a+1)v \relK u(a+1)wav
\end{equation}
for all $a\in\haken{n-1}$, $u,v,w\in\N^*$ such that $uaw(a+1)v\in
\Sym_n$ and $w$ contains the letter $a-1$ or the letter $a+2$;
\begin{equation}
  \label{rell}
  uav
  \relL
  (-1)^{\length{u}}
  a(\ol{u}\shuffle v)
\end{equation}
for all $a\in\haken{n}$, $u,v\in\N^*$ such that $uav\in \Sym_n$.
For later use, note that applying \kref{rell} twice gives 
\begin{equation}
  \label{x}
  aubv
  \relL
  (-1)^{n-1}
  \ol{v}b\ol{u}a
  \relL
  (-1)^{n-1+\length{v}}
  b(v\shuffle \ol{u}a)
\end{equation}
whenever $a,b\in\haken{n}$ and $u,v\in\N^*$ such that $aubv\in \Sym_n$.

\begin{remark}
  The space $\Lie_n^\perp$ is linearly generated by all non-trivial
  shuffles $u\shuffle v$, where $u,v\in\N^*$ such that $uv\in\Sym_n$
  (see, for instance, \cite{duchamp91}).
  As a consequence of Corollary~\ref{lie-char}, for fixed $a\in\haken{n}$,
  the elements
  $$
  uav-(-1)^{|u|}a(\ol{u}\shuffle v),
  $$
  where $u,v\in\N^*$ such that $uav\in\Sym_n$ and $u\neq\emptyword$,
  constitute a linear basis of $\Lie_n^\perp$. Another basis has been
  introduced by Duchamp (ibid.). This was pointed out to me by
  Christophe Reutenauer.
\end{remark}

This section concludes with a helpful observation concerning the order
reversing involution $\varrho_n=n(n-1)\ldots 1\in \Sym_n$.

\begin{prop} \label{inverter}
  $\varrho_n\Lie_n+\Lie_n\varrho_n\subseteq\Lie_n$,
  and
  $\varrho_n\QQ_n+\QQ_n\varrho_n\subseteq\QQ_n$.
  
  In particular, $\pi\relLK\sigma$ implies
  $\pi\varrho_n\relLK\sigma\varrho_n$ and
  $\varrho_n\pi\relLK\varrho_n\sigma$, for all $\pi,\sigma\in \Sym_n$.
\end{prop}

\begin{proof}
  First, $\varrho_n\omega_n=(-1)^{n-1}\omega_n$ yields
  $\varrho_n\Lie_n\subseteq\Lie_n$, while
  $\Lie_n\varrho_n\subseteq\Lie_n$ is obvious; and second, if
  $\sigma,\pi\in\Sym_n$ such that $\sigma\relK\pi$, then
  $\sigma\varrho_n\relK \pi\varrho_n$ and
  $\varrho_n\sigma\relK\varrho_n\pi$, as is readily seen from
  Proposition~\ref{coplactic-char}.  This implies
  $\varrho_n\QQ_n\subseteq\QQ_n$ and $\QQ_n\varrho_n\subseteq\QQ_n$.
  
  In particular, it follows that 
  $
  \varrho_n(\Lie_n\cap\QQ_n)^\perp+(\Lie_n\cap\QQ_n)^\perp\varrho_n 
  \subseteq
  (\Lie_n\cap\QQ_n)^\perp
  $,
  since $\varrho_n$ is an involution.
\end{proof}

%
%

\section{Coplactic Lie Elements} \label{copl-lie}

The aim of this section is to show $ \Lie\cap\QQ\subseteq\DD $, which
implies Theorem~\ref{intersect}, as was mentioned in the previous
section.  Throughout, $n\in\N$ is fixed.  Bearing in mind
Proposition~\ref{descent-char}, it suffices to show that
\begin{equation}
  \label{suffices}
  uaw(a+1)s
  \relLK
  u(a+1)was
\end{equation}
whenever $a\in\haken{n-1}$, $u,v,w\in\N^*$ such that $uaw(a+1)v\in
\Sym_n$ and $w\neq\emptyword$.

The essential idea of the proof is illustrated by the following

\begin{example}
  Let $\pi=15234,\sigma=25134\in\Sym_5$, then $\pi$ and $\sigma$ are
  in descent, but not in coplactic relation.  Applying \kref{rell}
  yields, however,
  \begin{eqnarray*}
    \pi
    & \relL &
    -5(1\shuffle 234)\\
    & = &
    -51234-52134-52314-52341\\
    & \relK &
    -51234-52134-51324-51342\\
    & = &
    -5(2\shuffle 134)\\
    & \relL &
    \sigma,
  \end{eqnarray*}
  hence $\pi\relLK\sigma$.
\end{example}

Some additional preparations are needed for the proof of
\kref{suffices}.  $v\in\N^*$ is called a \emph{sub-word} of
$w=w_1\ldots w_m\in\N^*$ if there exist $k\in\haken{m}$ and $1\le
i_1<\cdots<i_k\le m$ such that $ v=w_{i_1}\ldots w_{i_k} $.  For
instance, $23$ is a sub-word of $52143$.

For $a\in\haken{n-1}$, denote by $\tau_a=(a\;a+1)$ the transposition
in $\Sym_n$ swapping $a$ and $a+1$.  The word $v$ \emph{allows the
  $a$-switch in $\Sym_n$} if $\pi\relLK\pi\tau_a$ for all
$\pi=uaw(a+1)s\in \Sym_n$ such that $v$ is a sub-word of $w$. For
instance, $v=a+2$ and $v=a-1$ allow the $a$-switch in $\Sym_n$, by
\kref{relq}.  To save trouble, let it be said that, if $v$ contains a
letter twice or a letter $b>n$ or $b\in\{a,a+1\}$, then $v$ allows the
$a$-switch in $\Sym_n$; for in this case, there is no permutation
$\pi=uaw(a+1)s\in \Sym_n$ such that $v$ is a sub-word of $w$.  Another
way of stating \kref{suffices} now is that $v\in\N^*$ allows the
$a$-switch in $\Sym_n$ whenever $v\neq\emptyword$. The following three
helpful observations will be applied frequently.

\begin{prop} \label{beginner}
  Let $v\in\N^*$ such that $ \pi\relLK\pi\tau_a $ for all
  $$
  \pi=aw(a+1)s\in \Sym_n
  $$
  such that $v$ is a sub-word of $w$, then $v$ allows the $a$-switch
  in $\Sym_n$.
\end{prop}

\begin{proof}
  Let $\pi=uaw(a+1)s\in \Sym_n$ such that $v$ is a sub-word of $w$, then
  $$
  \pi
  \relL
  (-1)^{\length{u}}a(\ol{u}\shuffle w(a+1)s),
  $$
  by \kref{rell}.
  Each summand in this shuffle product is of the form
  $a\hat{w}(a+1)\hat{s}$ such that $w$ (hence also $v$) is a 
  sub-word of $\hat{w}$. It follows that
  $$
  \pi
  \relLK
  (-1)^{\length{u}}(a+1)(\ol{u}\shuffle was)
  \relL
  u(a+1)was,
  $$
  hence $v$ allows the $a$-switch in $\Sym_n$.
\end{proof}

\begin{prop} \label{backwards}
  Let $v=v_1\ldots v_m\in\N^*$ and assume that $v$ allows the
  $a$-switch in $\Sym_n$, then so does $\ol{v}$.
  Furthermore, if
  $$
  \tilde{v}
  :=
  (n+1-v_1)\ldots(n+1-v_m)
  \in
  \N^*,
  $$
  then $\tilde{v}$ allows the $(n-a)$-switch in $\Sym_n$.
\end{prop}

This is an immediate consequence of Proposition~\ref{inverter}.

\begin{prop} \label{general}
  Let $\pi\in \Sym_n$ and $\varphi_0,\ldots,\varphi_m\in \Z\Sym_n$ such that
  \begin{listing}
    \item 
      $\varphi_0=\pi$,
    \item
      $\varphi_i\relL \varphi_{i+1}$, or
      $\varphi_i\relK \varphi_{i+1}$ 
      and 
      $\varphi_i\tau_a\relK \varphi_{i+1}\tau_a$, 
      for all $i\in\haken{m-1}\cup\{0\}$,
    \item
      $\varphi_m\relLK\varphi_m\tau_a$,
  \end{listing}
  then $\pi\relLK\pi\tau_a$.
\end{prop}

\begin{proof}
  $\varphi\relL\psi$ implies $\varphi\tau_a\relL\psi\tau_a$ for all
  $\varphi,\psi\in \Z\Sym_n$, since $\Lie_n\tau_a=\Lie_n$.
  Combined with (ii), this implies 
  $\varphi_i\tau_a\relLK \varphi_{i+1}\tau_a$ for all
  $i\in\haken{m-1}\cup\{0\}$, hence
  $$
  \pi
  =
  \varphi_0
  \relLK
  \varphi_1
  \cdots
  \relLK
  \varphi_m
  \relLK
  \varphi_m\tau_a
  \relLK
  \cdots
  \relLK
  \varphi_1\tau_a
  \relLK
  \varphi_0\tau_a
  =
  \pi\tau_a\,,
  $$
  by (i) and (iii).
\end{proof}

We now show in four steps that each $v\in\N^*\setminus\{\emptyword\}$
allows the $a$-switch in $\Sym_n$.

\begin{step} \label{1}
  Let $a,k\in\haken{n-1}$ such that $k>a+1$, then $k(k+1)$ and
  $(k+1)k$ allow the $a$-switch in $\Sym_n$.
\end{step}

\begin{proof}
  Let $\pi = aw(a+1)s \in \Sym_n$ such that $k(k+1)$ or $(k+1)k$ is a
  sub-word of $w$.  It suffices to prove $\pi\relK \pi\tau_a$, by
  Proposition~\ref{beginner}.

  If $k=a+2$, then this follows from \kref{relq}.
  Let $k\ge a+3$, and proceed by induction on $k$.
  
  If $k-1$ occurs in $w$, then $\pi\relLK \pi\tau_a$, by induction.
  Let $k-1$ occur in $s$.
  
  If $k=a+3$ and $\pi=au_1(a+3)u_2(a+4)u_3(a+1)u_4(a+2)u_5$,
  then
  $$
  \pi
  \relK 
  a\,u_1\,(a+2)\,u_2\,(a+4)\,u_3\,(a+1)\,u_4\,(a+3)\,u_5\,.
  $$
  Applying Proposition~\ref{general}, yields
  $\pi\relLK\pi\tau_a$ in this case.
  In particular, $(a+3)(a+4)$ allows the $a$-switch in $\Sym_n$,
  hence also $(a+4)(a+3)$, by Proposition~\ref{backwards}.
  
  Now let $k>a+3$.  If $k-2$ occurs in $w$, then there are
  $u_i\in\N^*$ ($i\in\haken{5}$) such that either
  \begin{eqnarray*}
    \pi
    & = &
    au_1ku_2(k-2)u_3(a+1)u_4(k-1)u_5\\[1mm]
    & \relK &
    au_1(k-1)u_2(k-2)u_3(a+1)u_4ku_5=:\varphi
  \end{eqnarray*}
  or
  \begin{eqnarray*}
    \pi
    & = &
    au_1(k-2)u_2ku_3(a+1)u_4(k-1)u_5\\[1mm]
    & \relK &
    au_1(k-1)u_2ku_3(a+1)u_4(k-2)u_5=:\varphi\,.
  \end{eqnarray*}
  In both cases, $\varphi\relLK\varphi\tau_a$, by induction,
  hence $\pi\relLK\pi\tau_a$, by Proposition~\ref{general}.
  
  Assume that $k-2$ occurs in $s$, then there are $u_i\in\N^*$
  ($i\in\haken{6}$) such that one of the following four cases holds.

  \textbf{case 1.} 
  $\pi=au_1(k+1)u_2ku_3(a+1)u_4(k-2)u_5(k-1)u_6$, then
  $\pi\relK\pi\tau_{k-1}\relK\pi\tau_{k-1}\tau_k$,
  and 
  $$
  \pi\tau_{k-1}\tau_k
  =
  au_1ku_2(k-1)u_3(a+1)u_4(k-2)u_5(k+1)u_6
  \relLK
  (\pi\tau_{k-1}\tau_k)\tau_a\,,
  $$
  by induction. Again
  Proposition~\ref{general} implies $\pi\relLK\pi\tau_a$.

  \textbf{case 2.} 
  $\pi=au_1ku_2(k+1)u_3(a+1)u_4(k-1)u_5(k-2)u_6$, then
  put $m:=\length{u_4}+\length{u_5}+\length{u_6}+2$
  and 
  $$
  \varphi
  :=
  (-1)^{n-1+m}
  (a+1)\Big(
  u_4ku_5(k-2)u_6
  \shuffle 
  \ol{u}_3(k+1)\ol{u}_2(k-1)\ol{u}_1a
  \Big)
  $$
  to obtain
  $
  \pi 
  \relK 
  au_1(k-1)u_2(k+1)u_3(a+1)u_4ku_5(k-2)u_6
  \relL 
  \varphi
  $,
  by \kref{x}.  For all summands $\nu$ in $\varphi$ with
  $k$ to the left of $a$, that is $k\nu^{-1}<a\nu^{-1}$, $\nu\relLK
  \nu\tau_a$, by induction, while each of the summands $\nu$
  with $k$ to the
  right of $a$ is of the form
  \begin{eqnarray*}
    \nu
    & = &
    (a+1)v_1(k+1)v_2(k-1)v_3av_4kv_5
    \relK
    \nu\tau_k\,,
  \end{eqnarray*}
  hence $\nu\relLK\nu\tau_a$, by induction and
  Proposition~\ref{general}.  Putting both parts together yields
  $\varphi\relLK\varphi\tau_a$,
  hence $\pi\relLK\pi\tau_a$, by Proposition~\ref{general}.

  \textbf{case 3.} 
  $\pi=au_1(k+1)u_2ku_3(a+1)u_4(k-1)u_5(k-2)u_6$, then
  putting $m:=\length{u_4}+\length{u_5}+\length{u_6}+2$,
  \begin{eqnarray*}
    \pi
    & \relL &
    (-1)^{n-1+m}
    (a+1)\Big(
    u_4(k-1)u_5(k-2)u_6
    \shuffle 
    \ol{u}_3k\ol{u}_2(k+1)\ol{u}_1a
    \Big),
  \end{eqnarray*}
  by \kref{x}.  For each of the summands, swapping $a$ and
  $a+1$ yields an $\mathsf{LK}$-equivalent permutation, since either
  $k-1$ stands to the left of $a$ and the induction hypothesis may be
  applied, or $k-1$ stands to the right of $a$ and case 2 may be
  applied. Thus
  $\pi\relLK \pi\tau_a$, by Proposition~\ref{general}.
  
  Combining cases 1 and 3 shows that $(k+1)k$ allows the $a$-switch in
  $\Sym_n$, hence also $k(k+1)$, by Proposition~\ref{backwards}. This, in
  particular, yields the assertion in the remaining

  \textbf{case 4.} $\pi=au_1ku_2(k+1)u_3(a+1)u_4(k-2)u_5(k-1)u_6$.
\end{proof}

For the proof of the second step, an auxiliary result is needed.

\begin{prop} \label{Khelper}
  Let $j,k,x\in\haken{n}$, $v_1,\ldots,v_{k+1},w\in\N^*$ and set
  $x_j:=x+j$ for all $j\in\haken{k+1}$.
  If
  $$
  \pi
  :=
  v_1x_1v_2x_2
  \ldots 
  v_{j-1}x_{j-1}
  v_{k+1}x_{k+1}v_jx_j
  v_{j+1}x_{j+1}
  \ldots 
  v_kx_k\mbox{\hphantom{$_{+1}$}}
  w\hspace*{3pt}
  $$
  is contained in $\Sym_n$, then
  $$
  \pi
  \relK
  v_1x_2v_2x_3
  \ldots 
  v_{j-1}x_j\mbox{\hphantom{$_{-1}$}}
  v_{k+1}x_{j+1}v_jx_1
  v_{j+1}x_{j+2}
  \ldots 
  v_kx_{k+1}
  w\,.
  $$
\end{prop}

\begin{proof}
  Let
  $u:=v_1x_1v_2x_2\ldots v_{j-1}x_{j-1}$
  and 
  $\hat{u}:=v_{j+1}x_{j+2}\ldots v_kx_{k+1} w$, then
  \begin{eqnarray*}
    \pi
    & = &
    uv_{k+1}\ul{x_{k+1}}v_jx_jv_{j+1}x_{j+1}
    \ldots 
    v_{k-1}x_{k-1} v_k\ul{x_k}\mbox{\hphantom{$_{+1}$}} w\\[1mm]
    & \relK &
    u
    v_{k+1}\ul{x_k}\mbox{\hphantom{$_{+1}$}}v_jx_j v_{j+1}x_{j+1}
    \ldots 
    v_{k-1}\ul{x_{k-1}} v_kx_{k+1} w\\[1mm]
    & \relK &
    \cdots\\[1mm]
    & \relK &
    u v_{k+1}x_{j+1}v_jx_j v_{j+1}x_{j+2}\,
    \ldots 
    v_{k-1}x_k\mbox{\hphantom{$_{-1}$}} v_kx_{k+1} w\\[1mm]
    & = &
    v_1x_1v_2x_2
    \ldots 
    v_{j-2}x_{j-2}v_{j-1}\ul{x_{j-1}}v_{k+1}x_{j+1}
    v_j\ul{x_j}\mbox{\hphantom{$_{-1}$}}\hat{u}\\[1mm]
    & \relK &
    v_1x_1v_2x_2
    \ldots 
    v_{j-2}\ul{x_{j-2}}v_{j-1}x_j\mbox{\hphantom{$_{-1}$}}
    v_{k+1}x_{j+1}v_j\ul{x_{j-1}}\hat{u}\\[1mm]
    & \relK &
    \cdots\\[1mm]
    & \relK &
    v_1x_2v_2x_3
    \ldots 
    v_{j-2}x_{j-1}v_{j-1}x_j\mbox{\hphantom{$_{-1}$}}
    v_{k+1}x_{j+1}v_jx_1\mbox{\hphantom{$_{-1}$}}\hat{u}
  \end{eqnarray*}
  as asserted, where the letters in question are underlined in each
  step.
\end{proof}

\begin{step} \label{2}
  Let $a\in\haken{n-1}$ and $x,y\in\haken{n}$ such that $x,y>a+1$ or
  $x,y<a$, then $xy$ allows the $a$-switch in $\Sym_n$.
  
  In particular, if $v\in\N^*$ such that $\length{v}\ge 3$, then $v$ allows
  the $a$-switch in $\Sym_n$.
\end{step}

\begin{proof}
  If $x=y$, there is nothing to prove; let $x\neq y$.  Let
  $\pi=aw(a+1)s$ such that $xy$ is a sub-word of $w$.
  
  Consider first the case where $x,y>a+1$.  We may assume that $x<y$, by
  Proposition~\ref{backwards}.  
  The proof is done by induction on $m:=y-x$.

  If $m=1$, then $\pi\relLK\pi\tau_a$
  follows from Step~\ref{1}.
  
  Let $m>1$, and set $x_i:=x+i$ for all $i\in\N$.  Inductively, the
  case where $x_i$ occurs in $s$ for all $i\in\haken{m-1}$ remains.
  
  Choose $k\in\haken{m-1}$ maximal such that
  $
  x_1\pi^{-1}
  <
  x_2\pi^{-1}
  <
  \cdots
  <
  x_k\pi^{-1}
  $,
  that is
  $$
  \pi
  =
  a\,u_1\,x\,u_2\,y\,u_3\,(a+1)\,v_1\,x_1\,v_2\,x_2\,\ldots\,v_k\,x_k\,v_{k+1}
  $$
  for suitably chosen $u_i,v_i\in\N^*$.
  
  If $x_k=y-1$, then either $m=2$ and $\pi \relK \pi\tau_x$, or $m>2$
  and $\pi \relK \pi\tau_{y-1}$; in both cases, $\pi\relLK\pi\tau_a$,
  by induction and Proposition~\ref{general}.

  Let $x_k<y-1$, then $y>x_2$; and there is an index $j\in\haken{k}$ such that
  $x_{j-1}\pi^{-1}<x_{k+1}\pi^{-1}<x_j\pi^{-1}$ (where $x_0:=x$ if
  $j=1$).  Let $t:=au_1xu_2yu_3(a+1)$, then
  \begin{eqnarray*}
    \pi
    & = &
    tv_1x_1v_2x_2
    \ldots 
    v_{j-1}x_{j-1}
    v_j^{(1)}x_{k+1}v_j^{(2)}x_j
    v_{j+1}x_{j+1}
    \ldots 
    v_kx_k
    v_{k+1}\\[1mm]
    & \relK &
    tv_1x_2v_2x_3
    \ldots 
    v_{j-1}x_j
    v_j^{(1)}x_{j+1}v_j^{(2)}x_1
    v_{j+1}x_{j+2}
    \ldots 
    v_kx_{k+1}
    v_{k+1}\\[1mm]
    & =: &
    \hat{\pi}\,,
  \end{eqnarray*}
  by Proposition~\ref{Khelper}.
  Furthermore, $\hat{\pi}\relK \hat{\pi}\tau_x$, since
  $x\hat{\pi}^{-1}<x_2\hat{\pi}^{-1}<x_1\hat{\pi}^{-1}$.
  By induction,
  $\hat{\pi}\tau_x\relLK (\hat{\pi}\tau_x)\tau_a$, 
  hence also
  $\hat{\pi}\relLK \hat{\pi}\tau_a$, by Proposition~\ref{general}.
  Another application of Proposition~\ref{general} yields
  $\pi\relLK \pi\tau_a$
  and completes the proof in the case $x,y>a+1$.
  
  Now assume that $x,y<a$, then $n+1-x,n+1-y>(n-a)+1$, hence
  $(n+1-x)(n+1-y)$ allows the $(n-a)$-switch in $\Sym_n$, by the part
  already proven.  As a consequence, $xy$ allows the $a$-switch in
  $\Sym_n$, by Proposition~\ref{backwards}.
  
  If $v\in\N^*$ such that there are three distinct letters $\neq
  a,a+1$ occurring in $v$, then at least two of these are $<a$ or
  $>a+1$.  This completes the proof.
\end{proof}

\begin{step} \label{3}
  Let $a\in\haken{n-1}$ and $v\in\N^*$ such that $\length{v}\ge 2$, then $v$
  allows the $a$-switch in $\Sym_n$.
\end{step}

\begin{proof}
  By \kref{relq} and Step~\ref{2}, the case where $v=xy$ such that 
  and $x<a-1$ and $y>a+2$, or $y<a-1$ and $x>a+2$, remains.
  By Proposition~\ref{backwards}, it suffices to consider the case of
  $x<a-1$ and $y>a+2$.
  Let
  $
  \pi=aw(a+1)s\in \Sym_n
  $
  such that $xy$ is a sub-word of $w$.
  If $\length{w}\ge 3$, then $\pi\relLK\pi\tau_a$ follows from Step~\ref{2}.
  
  Let $w=xy$, then each of the letters $x+1,x+2,\ldots,a-1$ occurs in
  $s$, and $y>a+1>a>a-1\ge x+1$.  Applying Step~\ref{2} a number of
  times implies
  $$
  \pi
  \relLK
  \pi\tau_x
  \relLK
  \pi\tau_x\tau_{x+1}
  \relLK
  \cdots
  \relLK
  \pi\tau_x\tau_{x+1}\cdots\tau_{a-2}
  =
  a(a-1)y(a+1)\hat{s}
  $$
  for a properly chosen $\hat{s}\in\N^*$, hence $\pi\relLK\pi\tau_a$
  as asserted, by Proposition~\ref{general}.
\end{proof}

As the final step, we are now in a position to state and prove

\begin{thm} \label{copl-lie-thm}
  $\Lie\cap\QQ\subseteq \DD$.
\end{thm}

\begin{proof}
  It suffices to prove $\Lie_n\cap\QQ_n\subseteq\DD_n$, since
  $\Lie\cap\QQ=\bigoplus_{n\ge 0} \Lie_n\cap\QQ_n$.  By
  Proposition~\ref{beginner}, Step~\ref{3} and \kref{suffices}, it
  thus remains to be shown that $\pi\relLK\pi\tau_a$ whenever
  $x,a\in\haken{n}$, $s\in\N^*$ such that $\pi=ax(a+1)s\in \Sym_n$.
  
  If $n=3$, this is immediate. Let $n>3$ and choose $y\in\N$ and
  $u\in\N^*$ such that $s=yu$, then
  \begin{eqnarray*}
    \pi
    & = &
    ax(a+1)yu\\
    & = &
    a((a+1)\shuffle xyu)
    -
    a(a+1)xyu
    -
    axy((a+1)\shuffle u)\\
    & \relL &
    -
    (a+1)axyu
    -
    a(a+1)xyu
    -
    axy((a+1)\shuffle u)
    \quad\mbox{, by \kref{rell}}\\
    & \relLK &
    -
    (a+1)axyu
    -
    a(a+1)xyu
    -
    (a+1)xy(a\shuffle u)
    \quad\mbox{, by Step~\ref{3}}\\
    & \relL &
    (a+1)(a\shuffle xyu)
    -
    (a+1)axyu
    -
    (a+1)xy(a\shuffle u)
    \quad\mbox{, by \kref{rell}}\\
    & = &
    (a+1)xayu\\
    & = &
    \pi\tau_a\,.
  \end{eqnarray*}
  The theorem is proved.
\end{proof}

Denote by $\Delta$ the coproduct in $\PP$ (as in
\cite[p.~977]{malvenuto-reutenauer95}), and by $\emptyset$ the unique
element in $\Sym_0$.  For any a Hopf sub-algebra $\AA$ of $\PP$, let
$$
\Prim(\AA)
=
\setofall{\alpha\in\AA}{%
  \Delta(\alpha)=\alpha\tensor\emptyset+\emptyset\tensor \alpha}
$$
be the primitive Lie algebra of $\AA$, and denote by $\Prim(\AA)_n$ its
$n$-th homogeneous component.  The sub-algebra $\CC$ of $\PP$
generated by $\Prim(\PP)$ (the domain of co-commutativity of $\PP$)
contains $\LL$.  Furthermore, $\CC\cap \QQ$ (the domain of
co-commutativity of $\QQ$) is generated by $\Prim(\QQ)$ and contains
$\DD=\LL\cap \QQ$.

It turns out that $\DD$ is \emph{strictly} contained in $\CC\cap\QQ$.
Indeed, $\varphi=3412+2413-(3142+2143)\in\QQ_4\setminus\DD_4$, and
$\Delta(\varphi)=\varphi\tensor\emptyset+\emptyset\tensor \varphi$.
For $n=4$, $5$, $6$, the dimension of $\Prim(\DD)_n$ is, respectively,
$3$, $6$, and $9$, while the dimension of $\Prim(\QQ)_n$ is,
respectively, $4$, $9$, and $26$.  A description of the elements of
$\Prim(\DD)_n$ as well as of its dimension is known in general
(\cite[4.5]{blelau93a}, \cite[1.5]{blelau96}).  It would be of
interest if analogous results for $\QQ$ were obtained.

\newcommand{\etalchar}[1]{$^{#1}$}
 \newcommand{\noopsort}[1]{} \newcommand{\schocker}{ (= M. Schocker)}
  \newcommand{\appears}[1]{to appear in \emph{#1}}

\end{document}